\documentclass[12pt,a4paper]{article}
\usepackage{amssymb,amsmath}
\usepackage{graphicx}	

\newlength{\IndentI}
\newlength{\IndentII}
\newlength{\IndentIII}
\newlength{\WidthI}
\newlength{\WidthII}
\newlength{\WidthIII}
\newcommand{\eventone}{E_1}
\newcommand{\eventtwo}{E_2}
\newcommand{\eventthree}{E_3}
\newcommand{\eventfour}{E_4}
\newcommand{\eventfive}{E_5}
\newcommand{\eventsix}{E_6}

\def\P{\mathcal{P}}
\def\Q{\mathcal{Q}}
\def\K{\mathcal{K}}
\def\ulrmor{\underline{\mathrm{or}}}

\newtheorem{theorem}{Theorem}

\newtheorem{lemma}{Lemma}
\newtheorem{corollary}{Corollary}
\newenvironment{proof}
{\par\noindent\textit{{\bfseries Proof:}}~ } 
{\hfill\rule{2.1mm}{2.1mm}\bigskip\par} 

\begin{document}

\title{Implications of contrarian and one-sided strategies for the
fair-coin game}
\author{Yasunori Horikoshi and Akimichi Takemura\\ \\
        Department of Mathematical Informatics\\
        Graduate School of Information Science and Technology\\
        The University of Tokyo\\
        }
\date{March, 2007}

\maketitle

\begin{abstract}
  We derive some results on contrarian and one-sided strategies by
  Skeptic for the fair-coin game in the framework of the game-theoretic
  probability of Shafer and Vovk \cite{sv}.  In particular, concerning
  the rate of convergence of the strong law of large numbers (SLLN),
  we prove that Skeptic can force that the convergence has to be
  slower than or equal to $O(n^{-1/2})$.  This is achieved by a very
  simple contrarian strategy of Skeptic.  
  This type of result, bounding the rate of convergence from
  below, contrasts with more standard results of bounding the rate of
  SLLN from above by using momentum strategies.  We also derive a
  corresponding one-sided result.
\end{abstract}

\section{Introduction}
\label{sec:intro}

In the theory of game-theoretic probability by Shafer and Vovk
\cite{sv}, various ``probability laws'' such as SLLN are proved by
constructing explicit strategies of Skeptic, who is one of the two
players in a game.  Construction of a clever and explicit strategy of
Skeptic often leads to a remarkably simple proof of the corresponding
result in the measure theoretic probability theory, even without
preparations from measure theory.  This is already apparent in the
simple strategy  used in Chapter 3 of Shafer and Vovk \cite{sv}, where
Skeptic always bets a fixed proportion of his capital.
See also the Bayesian strategies of Skeptic in coin-tossing games 
in \cite{ktt}.
New problems and
their solutions offered by the framework of the game-theoretic
probability are now actively investigated in various directions.
For background material and further developments of the game-theoretic
probability see Vovk and Shafer \cite{vovk/shafer:2005} and Shafer and
Vovk \cite{source1} and references therein.
See Takeuchi \cite{takeuchi:2004} for some original ideas and results.
Defensive forecasting, which is a new non-parametric forecasting 
method based on the game-theoretic probability, was initiated in 
\cite{vovk-takemura-shafer-defensive} and 
\cite{wp10}.

We can roughly classify strategies of Skeptic into two classes, namely
the class of momentum strategies and the class of contrarian
strategies.  This distinction was clearly demonstrated 
in  a talk by Glenn Shafer \cite{shafer-at-u-tokyo}.  See also
\cite{root-dt}.
Consider again the convergence in SLLN.
In momentum strategies Skeptic assumes that Reality, the other player
of the game,  will
keep deviating from the (zero) theoretical mean
in the same direction and bets
accordingly. In contrast, in contrarian strategies Skeptic assumes
that Reality tries to decrease deviation from the mean and bets
accordingly.  In both of these strategies, Skeptic looks only at the
absolute deviation and in this sense these strategies
are two-sided strategies.  A more primitive strategy of Skeptic is
one-sided and bets only toward a particular direction (up or down).
In Shafer and Vovk \cite{sv} one-sided strategies are treated as
restrictions on the move space of Skeptic, namely Skeptic is only
allowed to buy a certain ``ticket''.  However as we discuss in Section
\ref{sec:one-sided}, Skeptic can use one-sided strategies to force
stronger unbiasedness to Reality than implied by two-sided
strategies.

It is only natural to expect that stronger results require more
complicated strategies by Skeptic.  
For example in \cite{kt}  we have shown that 
a simple strategy of Skeptic based on the past average of the moves by
Reality forces SLLN 
for the case of bounded Reality's moves.  However 
if Reality's moves are unbounded, strategies for forcing SLLN are much
more complicated as discussed in \cite{ktt-unbounded}.

As another example, we mention that the proof of the law of the iterated
logarithm (LIL) in Chapter 5 of Shaver and Vovk \cite{sv} is much
harder than the proof of SLLN or  the central limit theorem.
There are two parts in the proof of LIL.  
In the first part the growth rate of the capital process
is bounded from above by a momentum strategy 
and in the second part it is bounded from below by a contrarian strategy.
As seen from the proof, the latter part
of the proof is much more difficult.  This suggests that construction
of effective contrarian strategies of Skeptic is more challenging than
construction of momentum strategies.  In this paper we only consider
the fair-coin game, which is the simplest game in the game-theoretic
probability.  However we believe that the contrarian
strategies obtained in this paper can be generalized further and give insights
on other strategies in various more general games in the game-theoretic
probability. We should also mention that our results on the rate of
convergence of SLLN are already implied by LIL.  Therefore the merit
of this paper is to clarify implications of very simple explicit contrarian
strategies.

Organization of the paper is as follows.
For the rest of this section, we briefly introduce necessary notations
and definitions from \cite{sv}. 
In Section \ref{sec:contrarian} we consider contrarian strategies.
It is divided into two subsections.  In subsection
\ref{sec:multiplicative-contrarian} multiplicative contrarian
strategies based on the past average of Reality's moves are studied 
and in subsection
\ref{sec:additive-contrarian} additive
contrarian strategies based on the past sum of Reality's moves are 
studied.  In Section \ref{sec:one-sided} we study one-sided strategies
to strengthen results obtained in Section \ref{sec:contrarian}. We end
the paper with some discussions in 
Section \ref{sec:discussion}.

\subsection{Notations and definitions}
Here we summarize necessary notations
and definitions from \cite{sv} for our paper.
We also give a 
formal definition of stopping times.

In this paper we consider the fair-coin game.  Its protocol is given
as follows.
\bigskip

\parshape=7
\IndentI   \WidthI
\IndentI   \WidthI
\IndentI   \WidthI
\IndentII  \WidthII
\IndentII  \WidthII
\IndentII  \WidthII
\IndentI   \WidthI
\noindent
{\scshape Fair-Coin Game}\\
$\K_0 :=1$.\\
FOR  $n=1, 2, \dots$\,:\\
  Skeptic announces $M_n\in\mathbb{R}$.\\
  Reality announces $x_n\in\{-1, 1\}$.\\
  $\K_n := \K_{n-1} + M_n x_n$.\\
END FOR

\medskip
A finite sequence $t = x_1x_2\cdots x_n$ 
consisting of $1$ and $-1$ is called a situation and the
set of all situations is denoted by $\Omega^\diamondsuit$. 
The length of $t = x_1x_2\cdots x_n$ is $n$.
The initial situation, which is the special situation of the
length 0, is
denoted by $\Box \in \Omega^\diamondsuit$.  If a situation $t$
is an initial segment of another situation $t'$, we say that $t$ precedes
$t'$ and $t'$ follows $t$. For a situation $t=x_1x_2\cdots x_n$, we define
$-t:=(-x_1)(-x_2)\cdots(-x_n)$.
If we consider the infinite binary tree describing the progress
of the fair-coin game, a situation can be identified with a node in the
tree.
 
An infinite
sequence $\xi=x_1 x_2 \dots$ 
consisting of 1 and $-1$ is called a path and the set of all
paths is denoted by $\Omega$. An event $E$ is a subset of $\Omega$.
For a path $\xi$, the situation consisting of the first $n$ terms of
$\xi$ is written as $\xi_n$. 
Given a path $\xi$ and a situation $t$, if there exists $n$ such that
 $\xi_n=t$, we say that $\xi$ goes through $t$. Given a
 situation $t$, we define the cylinder set $O_t\subset \Omega$ by
 \[
 O_t:=\{\xi \mid \xi\ \mathrm{goes\ through}\ t\}.
 \]

A process is a 
function $\Omega^\diamondsuit\rightarrow
\mathbb{R}$ and  a variable is 
a function $\Omega\rightarrow\mathbb{R}$.
Given a process $f$, a variable $f_n$ is defined by
\begin{equation}
\label{eq:process-variable}
 f_n(\xi):=f(\xi_n).
\end{equation}
In this paper the symbols $s$ and $\overline{x}$ are used
for two special processes, the sum and the average. They are defined by
\begin{eqnarray*}
 s(x_1x_2\cdots x_n) &:=& x_1+x_2+\cdots + x_n,\\
 \overline{x}(x_1x_2\cdots x_n) &:=& \frac{x_1+x_2+\cdots +x_n}{n},
\end{eqnarray*}
where $s(\Box)=\overline{x}(\Box)=0$.  By 
(\ref{eq:process-variable}), we also write
$s(x_1x_2\cdots x_n)=s_n(\xi)$, 
$\overline{x}(x_1x_2\cdots x_n) = \overline{x}_n(\xi)$.

The strategy $\P$ of Skeptic is a process. 
When $\P$ is a strategy of Skeptic, the capital
process of $\P$ (with zero initial capital) is 
denoted by $\K^\P$ and defined by 
\[
 \K^\P(x_1x_2\cdots x_n) := 
  \P(\Box)x_1+\P(x_1)x_2+\cdots +\P(x_1x_2\cdots x_{n-1})x_n
\]
and $\K^\P(\Box):=0$.  
The capital process of $\P$ is the total amount Skeptic earns when he
follows the strategy $\P$. 
If Skeptic uses the strategy $\P$ with the initial
capital $a$ and Reality chooses the path $\xi$, then the capital Skeptic
holds at the end of the $n$-th round is $a+\K_n^\P(\xi)$. When $\xi$ is
fixed, we write simply $s_n, \overline{x}_n$ or $\K_n^\P$.


We say that Skeptic can weakly force an event $E\subset\Omega$, if
there exists a strategy $\P$ of Skeptic such that
\begin{equation}
\label{eq:weakly-forcing}
\limsup_n \K_n^\P(\xi) = \infty, \qquad \forall \xi\not\in E, 
\end{equation}
under the restriction of the ``collateral duty''
\[
\K_n^\P(\xi) \ge -1, \quad \forall \xi\in \Omega, \forall n\ge 0.
\]

We say that Skeptic can force $E$ if $\limsup_n$ in 
(\ref{eq:weakly-forcing}) is replaced by $\lim_n$.
By Lemma 3.1 of \cite{sv}, if Skeptic can weakly force $E$, then he can
force $E$.

The upper price $\overline{\mathrm{E}}[x]$ of a variable $x$ 
is defined by
\begin{equation}
 \overline{\mathrm{E}}[x] := \inf\{a \mid \exists \P \ 
\forall \xi \in \Omega: a+\K_n^{\P}(\xi) \geq x(\xi)\ \mathrm{a.a.}\},
\end{equation}
where a.a.\ (almost always) means ``except for a finite number of $n$''.
We can regard a variable $x$ as a ticket whose holder earns $x(\xi)$ if
Reality chooses $\xi$ and  $\overline{\mathrm{E}}[x]$ is the infimum of
the initial capital with which Skeptic can superreplicate $x$.
When $a \geq \overline{\mathrm{E}}[x]$, we say that
Skeptic can buy $x$ for $a$.
Given a situation $t$, the upper price of $x$ on the situation $t$ is also
defined and denoted by $\overline{\mathrm{E}}_t[x]$ (Chapter 1 of
\cite{sv}).

Finally we give a formal definition of a stopping time.
A stopping time is a variable $f$: $\Omega \rightarrow
\{1,2,\dots\}\cup \{\infty\}$,  which satisfies
\[
 \forall \xi, \xi'\in \Omega, \forall n\in \mathbb{N} : (f(\xi)=n\ \&\ \xi_n = \xi'_n) \Rightarrow f(\xi')=n,
\]
where $\mathbb{N}=\{1,2,\dots\}$.
For a stopping time $f$ and a path $\xi$, if $f(\xi)<\infty$,
$\xi_{f(\xi)}$ is the situation where the value of $f(\xi)$ is
determined.
We let $[f]$ denote 
the set of situations where the value of $f$ is determined:
\begin{equation}
\label{eq:finite-stopping-time}
 [f] := \{t\in\Omega^\diamondsuit \mid \exists \xi:\xi_{f(\xi)}=t\}.
\end{equation}
A stopping time $f$ can be identified with $[f]$.

\section{Contrarian strategies}
\label{sec:contrarian}

In this section we consider forcing the following events:
\begin{align}
 \eventone& := \{\xi  \mid
 \limsup_{n\rightarrow\infty} \sqrt{n}|\overline{x}_n| \geq 1   \},  \\ 
 \eventtwo &:= \{\xi \mid  |s_n| > \sqrt{n}-1\ \mathrm{i.o.} \},
\end{align}
where i.o.\ (infinitely often) means ``for infinitely many $n$''.
In Subsection 
\ref{sec:multiplicative-contrarian}
we prove that Skeptic can force $\eventone$ by
a mixture of multiplicative contrarian strategies
based on the past average of Reality's moves and 
in Subsection 
\ref{sec:additive-contrarian}
we prove that Skeptic can force
$\eventtwo$ by 
a mixture of additive
contrarian strategies based on the past sum of Reality's moves.
Since $\eventtwo\subset \eventone$, forcing $\eventtwo$ is stronger
than forcing $\eventone$.  However the multiplicative strategy in
Subsection
\ref{sec:multiplicative-contrarian}
is of interest, because it is a contrarian counterpart of the momentum 
strategy studied in \cite{kt}.

\subsection{Multiplicative contrarian strategy}
\label{sec:multiplicative-contrarian}

In this section we study the following multiplicative 
contrarian strategy $\P_c$:
\[
 \P_c : \ M_n = -c\overline{x}_{n-1}\K_{n-1},
\]
where $c$ is an arbitrary constant satisfying $0<c\leq \frac{1}{2}$
and the initial capital  is $\K_0=a=1$. The case of $-\frac{1}{2}
\leq c < 0$ was studied in \cite{kt}.
Let
\[
\Q=\sum_{i\in \mathbb{N}}  \frac{1}{2^i}\P_{1/2^i}
\]
denote an infinite mixture of the strategies $\P_c$, $c=1/2^i$,
$i=1,2,\dots$.  Then the following result holds.

\begin{theorem}\label{thm1}
 Skeptic can force $\eventone$ by $\Q$.
\end{theorem}

The rest of this subsection is devoted to a proof of Theorem \ref{thm1}.
Let
 \[
 \eventone^c := \left\{\xi\ \left|\ \limsup_{n\rightarrow\infty}(1+2c)\sqrt{n}|\overline{x}_n| \geq 1\right.\right\}. 
 \]
Then $\eventone$ can be represented as
\[
 \eventone = \bigcap_{i\in \mathbb{N}} \eventone^{1/2^i}.
\]
Thus, by Lemma 3.2 of \cite{sv},
Theorem \ref{thm1} is a consequence of the following lemma.
\begin{lemma}
\label{lem:for-thm-1}
Skeptic can force $\eventone^c$ with $\P_c$.
\end{lemma}


 In order to prove the lemma, we need to show that for any $\xi\in \Omega$
 \begin{itemize}
  \item Skeptic's capital $\K_n=1+\K_n^{\P_c}$ never gets negative, and
  \item $\limsup_{n\rightarrow\infty} (1+2c)\sqrt{n}|\overline{x}_n| <
	1\ \Rightarrow\ \K_n\rightarrow \infty \quad  (n\rightarrow \infty)$.
 \end{itemize}
 By definition for any $\xi$,
 \begin{eqnarray}
  \K_n &=& \K_{n-1}-c\overline{x}_{n-1}\K_{n-1}x_n\nonumber\\
  &=&\K_{n-1}(1-c\overline{x}_{n-1}x_n)\nonumber\\
  &=&\prod_{i=2}^n(1-c\overline{x}_{i-1}x_i).\label{Feb 09 14:39:30 2007}
 \end{eqnarray}
 In the expression (\ref{Feb 09 14:39:30 2007}) the index $i$ starts
 from 2 because $\overline{x}_0=0$. Also for any $i=2, 3, \ldots, n$,
 $1-c\overline{x}_{i-1}x_i>0$. Thus the first statement on the
 collateral duty is trivial.  

 \bigskip
 We divide the proof of the second statement into four parts. 

 \paragraph{\underline{Step 1}}
 In step 1 and step 2, we 
fix an arbitrary path $\xi = x_1x_2\cdots \in \Omega$.
 Since $\log (1+t)\geq t-t^2$ whenever $|t|\leq \frac{1}{2}$, from (\ref{Feb 09 14:39:30 2007})
 \begin{eqnarray}
 \log\K_n &=& \sum_{i=2}^{n}\log(1-c\overline{x}_{i-1}x_i)\nonumber\\
 &\geq&
  -c\sum_{i=2}^{n}\overline{x}_{i-1}x_i-c^2\sum_{i=2}^{n}\overline{x}_{i-1}^2x_i^2.\label{Aug 23 20:44:00 2006}  
 \end{eqnarray}
We use the identity
\begin{equation}
  \sum_{i=2}^{n}\overline{x}_{i-1}x_i =
 \frac{1}{2}\sum_{i=2}^{n}\frac{i}{i-1}\overline{x}_i^2
 +\frac{n}{2}\overline{x}_n^2-\frac{1}{2}\left(x_1^2+\sum_{i=2}^{n}\frac{1}{i-1}x_i^2\right)\label{Feb 09 14:56:13 2007}
\end{equation}
which is shown in \cite{kt} and easily follows from
\[
s_{n-1}x_n = \frac{1}{2} (s_n^2 - s_{n-1}^2 - x_n^2).
\]
Substituting (\ref{Feb 09 14:56:13 2007})
into (\ref{Aug 23 20:44:00 2006}), we have
\begin{align}
\log \K_n
  &\geq 
 -c\left\{\frac{1}{2}\sum_{i=2}^{n}\frac{i}{i-1}\overline{x}_i^2+\frac{n}{2}\overline{x}_n^2-\frac{1}{2}\left(x_1^2+\sum_{i=2}^{n}\frac{1}{i-1}x_i^2\right)\right\}-c^2\sum_{i=2}^{n}\overline{x}_{i-1}^2x_i^2\nonumber\\
 &=
  -\frac{c}{2}\sum_{i=2}^{n}\frac{i}{i-1}\overline{x}_i^2-\frac{nc}{2}\overline{x}_n^2+\frac{c}{2}\left(x_1^2+\sum_{i=2}^{n}\frac{1}{i-1}x_i^2\right)-c^2\sum_{i=2}^{n}\overline{x}_{i-1}^2x_i^2.\label{Sep 04 19:10:59 2006}
\end{align}
Now, $x_i^2=1$ because $x_i\in\{-1, 1\}$ and 
 \[
 1+\frac{1}{2}+\cdots+\frac{1}{n-1}\geq \log n =
 \int_{1}^{n}\frac{1}{x}\mathrm{d}x.
 \]
Thus we have
 \begin{eqnarray}
 (\ref{Sep 04 19:10:59 2006})&=&
 -\frac{c}{2}\sum_{i=2}^{n}\frac{i}{i-1}\overline{x}_i^2-\frac{nc}{2}\overline{x}_n^2+\frac{c}{2}\left(1+\sum_{i=2}^{n}\frac{1}{i-1}\right)-c^2\sum_{i=2}^{n}\overline{x}_{i-1}^2\nonumber\\
 &\geq&
  -\frac{c}{2}\sum_{i=2}^{n}\frac{i}{i-1}\overline{x}_i^2-\frac{nc}{2}\overline{x}_n^2+\frac{c}{2}\left(1+\log n\right)-c^2\sum_{i=2}^{n}\overline{x}_{i-1}^2\nonumber\\
  &=&
  \frac{c}{2}\left(1+\log n\right)-\left(\frac{c}{2}\sum_{i=2}^{n}\frac{i}{i-1}\overline{x}_i^2+c^2\sum_{i=2}^{n}\overline{x}_{i-1}^2+\frac{nc}{2}\overline{x}_n^2\right)\nonumber\\
  &=&\frac{c}{2}\left\{
  1+\log n-\left(\sum_{i=2}^{n}\frac{i}{i-1}\overline{x}_i^2+2c\sum_{i=2}^{n}\overline{x}_{i-1}^2+n\overline{x}_n^2\right)\right\}\nonumber\\
 &=&\frac{c}{2}\left\{
  1+\log n\left(1-\frac{\sum_{i=2}^{n}\frac{i}{i-1}\overline{x}_i^2+2c\sum_{i=2}^{n}\overline{x}_{i-1}^2+n\overline{x}_n^2}{\log n}\right)\right\}.\label{Oct 10 18:34:55 2006}
 \end{eqnarray}

\bigskip\noindent
By (\ref{Oct 10 18:34:55 2006}) we have shown that
\[
\xi \in F_c \ \Rightarrow \ 
\lim_{n \rightarrow \infty}\K_n=\infty.
\]
where

\begin{equation}
\label{Sep 03 14:00:11 2006}
F_c = \left\{ \xi \ \left|\ 
 \limsup_{n\rightarrow \infty}
    \frac{\sum_{i=2}^{n}\frac{i}{i-1}\overline{x}_i^2+2c\sum_{i=2}^{n}\overline{x}_{i-1}^2+n\overline{x}_n^2}{\log
    n}< 1 \right.\right\}.
\end{equation}

\paragraph{\underline{Step 2}}
Rewriting the numerator in (\ref{Sep 03 14:00:11 2006}), we have
 \begin{eqnarray}
 \hspace{1.5cm}&&\hspace{-2cm}\sum_{i=2}^{n}\frac{i}{i-1}\overline{x}_i^2+2c\sum_{i=2}^{n}\overline{x}_{i-1}^2+n\overline{x}_n^2\nonumber\\
 &=&2c\overline{x}_1^2+\sum_{i=2}^{n-1}\left(\frac{i}{i-1}+2c\right)\overline{x}_i^2+\frac{n}{n-1}\overline{x}_n^2+n\overline{x}_n^2\nonumber\\
 &=& 2c+\sum_{i=2}^{n-1}\left(1+\frac{1}{i-1}+2c\right)\overline{x}_i^2+\left(1+\frac{1}{n-1}+n\right)\overline{x}_n^2\nonumber\\
 &\leq&
  2c+\sum_{i=2}^{n-1}\left(1+\frac{1}{i-1}+2c\right)\overline{x}_i^2+n\overline{x}_n^2+\left(1+\frac{1}{n-1}\right),\label{Feb 09 15:39:03 2007} 
 \end{eqnarray}
where we used $\overline{x}_1^2=1$ and $\overline{x}_n\leq 1$.

\bigskip 
By (\ref{Feb 09 15:39:03 2007}) we can state that for any $\epsilon> 0$
there exist $N_1(\epsilon)$ and $A_\epsilon$ such that for any
$\xi \in \Omega$
\begin{equation}
   \sum_{i=2}^{n}\frac{i}{i-1}\overline{x}_i^2+2c\sum_{i=2}^{n}\overline{x}_{i-1}^2+n\overline{x}_n^2\leq(1+\epsilon+2c)\sum_{i=N_1(\epsilon)}^{n-1}\overline{x}_i^2+n\overline{x}_n^2+A_\epsilon.\label{Jan 20 15:57:55 2007}
\end{equation}

\paragraph{\underline{Step 3}}
Now we consider the following event. 
\[
F_{c,\epsilon}=
\{ \xi \mid 
  \limsup_{n\rightarrow \infty}\ (1+\epsilon+2c)n\overline{x}_n^2 < 1
\},
\]
where $\epsilon > 0$ is fixed.
We will show that  $F_{c,\epsilon}\subset F_c$.
Fix any path $\xi\in F_{c.\epsilon}$. 
Then there exist $\alpha = \alpha_\xi<1$ and $N_2 = N_2(\xi)$ such
that for any $n\geq N_2$
\[
 (1+\epsilon+2c)\overline{x}_n^2\leq \alpha\frac{1}{n}.
\]
Therefore for $n> N_2$
\begin{align}
(1+\epsilon+2c)\sum_{i=N'_1(\xi)}^{n-1}\overline{x}_i^2
  &=(1+\epsilon+2c)\sum_{i=N'_1(\xi)}^{N_2-1}\overline{x}_i^2+(1+\epsilon+2c)\sum_{i=N_2}^{n-1}\overline{x}_i^2\nonumber\\
&\leq
  (1+\epsilon+2c)\sum_{i=N'_1(\xi)}^{N_2-1}\overline{x}_i^2+\alpha\left(\frac{1}{N_2}+\cdots+\frac{1}{n-1}\right)\nonumber\\
 &\leq
  \alpha\log n+B_\xi,\label{Feb 09 16:50:01 2007}
\end{align}
where $N'_1(\xi) = N_1(\epsilon(\xi))$ and $B_\xi$ is a constant.
Substituting (\ref{Feb 09 16:50:01 2007}) into (\ref{Jan 20 15:57:55
2007}), we have
\[
\sum_{i=2}^{n}\frac{i}{i-1}\overline{x}_i^2+2c\sum_{i=2}^{n}\overline{x}_{i-1}^2+n\overline{x}_n^2\leq\alpha\log n+n\overline{x}_n^2+C_\xi,
 \]
 where $C_\xi = A_{\epsilon(\xi)}+B_\xi$. Thus 
\begin{align*}
\limsup_{n\rightarrow \infty}
    \frac{\sum_{i=2}^{n}\frac{i}{i-1}\overline{x}_i^2
  +2c\sum_{i=2}^{n}\overline{x}_{i-1}^2+n\overline{x}_n^2}{\log n}
& \leq\limsup_{n\rightarrow\infty}\frac{\alpha\log
   n+n\overline{x}_n^2+C_\xi}{\log n}\\
  &=\alpha  < 1,
\end{align*}
where we used $\displaystyle\limsup_{n\rightarrow
  \infty}\frac{n\overline{x}_n}{\log n} = 0$ since
  $\limsup_{n\rightarrow \infty}(1+\epsilon+2c)n\overline{x}_n^2<1$.

\paragraph{\underline{Step 4}}
By definition for any $\xi\in \Omega\backslash \eventone^c$ we have
 \[
 \limsup_{n\rightarrow \infty}\ (1+2c)n\overline{x}_n^2 < 1.
 \]
So we can find $\epsilon(\xi)$ such that
\[
 \limsup_{n\rightarrow \infty}\ (1+\epsilon(\xi)+2c)n\overline{x}_n^2 < 1.
\]
Then Step 1 and Step 3 show $\K_n\rightarrow \infty\ (n\rightarrow
\infty)$ if Reality chooses $\xi\in\Omega\backslash \eventone^c$. This
completes the proof of Lemma \ref{lem:for-thm-1}.


\subsection{Additive contrarian strategies}
\label{sec:additive-contrarian}
In this section we prove the following theorem which gives a somewhat stronger
statement than Theorem \ref{thm1}:

\begin{theorem}\label{thm2}
Skeptic can weakly force $\eventtwo$.
\end{theorem}

In order to prove this theorem we need to combine various strategies.
The basic ingredient is an additive contrarian strategy in
(\ref{Feb 26 16:27:11 2007}) below.  Other strategies will be studied
in separate subsections.
Here is the general flow of the proof. First (section \ref{sec1}) we
will construct a strategy which makes the initial capital $\K_0 = 1$
increase to $\K_n =
1+\frac{n}{2}\epsilon$ when $s_n=0$. 
Next (section \ref{sec2-onesided}), we will construct a strategy
weakly forcing $\{s_n = 0\ \mathrm{i.o.}\}$. These strategies have a
risk that the capital becomes 
negative if Reality makes $|s_n|$ as large as $\sqrt{n}$. Therefore 
Skeptic must stop running the strategies right before his capital becomes
negative in order to observe the collateral duty.
But if Reality keeps $|s_n|$ smaller than $\sqrt{n}$
forever, Skeptic
can keep running strategies and then $\limsup \K_n = \infty$. Thus
Skeptic can weakly force that $|s_n|$ becomes as large as $\sqrt{n}$
eventually. Now dividing the initial capital into countably many
accounts, Skeptic can weakly force that 
$|s_n|$ become as large as $\sqrt{n}$ infinitely often.

\subsubsection{The strategy increasing the capital when
    the sum process returns to the origin}
\label{sec1}
    
Here we consider the following additive 
contrarian strategy $\tilde{\P}^\epsilon$: 
\begin{equation}
   \tilde{\P}^\epsilon :\ 
   M_n = -\epsilon s_{n-1},\label{Feb 26 16:27:11 2007}
\end{equation}
where $\epsilon > 0$ is a small positive constant.
If we temporarily ignore the collateral duty, the 
strategy $\tilde{\P}^\epsilon$ has a very simple explicit 
capital process described in the next lemma.
\begin{lemma}\label{lem2}
\begin{equation}
   \K_n^{\tilde{\P}^\epsilon} = \frac{\epsilon}{2}(n-s_n^2).\label{Feb 25 22:53:19 2007}
\end{equation}
\end{lemma}

\begin{proof}
 We use an induction on $n$.
 When $n=0$, (\ref{Feb 25
 22:53:19 2007}) holds by definition. Now  assume that (\ref{Feb 25
 22:53:19 2007}) holds for $n=k$. There are five cases depending on the
 signs of $x_{k+1}$ and $s_k$.
 \begin{enumerate}
  \item $s_k = 0$.
  \item $x_{k+1}=1, s_k > 0$.
  \item $x_{k+1}=1, s_k < 0$.
  \item $x_{k+1}=-1, s_k > 0$.
  \item $x_{k+1}=-1, s_k < 0$.
 \end{enumerate}
 \paragraph{\boldmath case 1) $s_k=0$ :}
 By the assumption of induction,
 \[
  \K_{k}^{\tilde{\P}^\epsilon}=\frac{\epsilon}{2}k.
 \]
 Since $s_k=0$, $s_{k+1}^2=1$ and $M_{k+1}=0$. Thus we have
\begin{eqnarray*}
 \K_{k+1}^{\tilde{\P}^\epsilon} &=&\K_{k}^{\tilde{\P}^\epsilon}+M_{k+1}x_{k+1}\\
 &=& \frac{\epsilon}{2}k \\
 &=&\frac{\epsilon}{2}(k+1-s_{k+1}^2).
\end{eqnarray*}

 \paragraph{\boldmath case 2) $x_{k+1} = 1, s_k>0$ :}
 By the assumption of induction,
 \[
   \K_{k}^{\tilde{\P}^\epsilon}=\frac{\epsilon}{2}(k-s_k^2).
 \]
 Since $M_{k+1}=-\epsilon s_k$ and $s_{k+1}^2=(s_k+1)^2$, we have
 \begin{eqnarray*}
  \K_{k+1}^{\tilde{\P}^\epsilon} &=&\K_{k}^{\tilde{\P}^\epsilon}+M_{k+1}x_{k+1}\\
  &=& \frac{\epsilon}{2}(k-s_k^2)-\epsilon s_k\\
  &=& \frac{\epsilon}{2}(k-s_k^2-2 s_k)\\
  &=&\frac{\epsilon}{2}(k+1-(s_k^2+2 s_k+1))\\
  &=&\frac{\epsilon}{2}(k+1-(s_k+1)^2)\\
  &=&\frac{\epsilon}{2}(k+1-s_{k+1}^2).
 \end{eqnarray*}
The other cases are proved by almost the same argument.
\end{proof}

\begin{figure}[h]
 \hspace{-1.5cm}
 \begin{center}
  \includegraphics[height=0.6\textwidth]{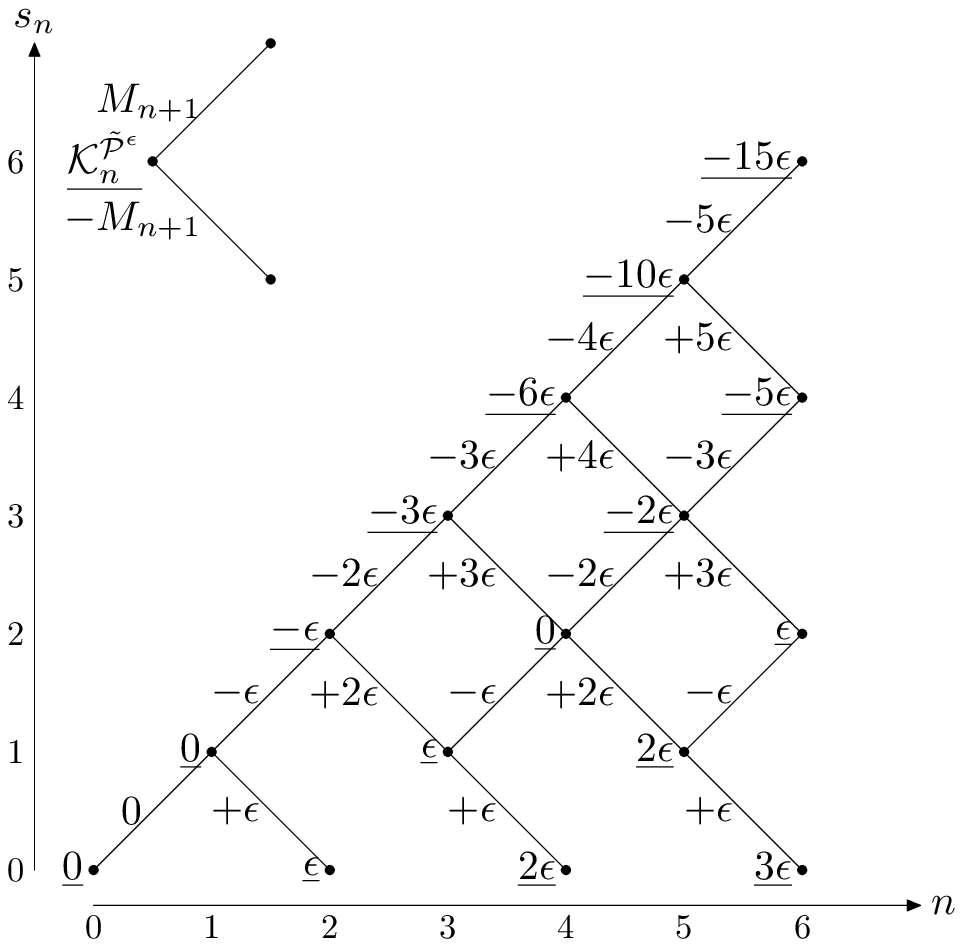}
 \end{center}
 \caption{Behavior of $\K_n^{\tilde{\P}^\epsilon}$} \label{fig}
\end{figure}

We can intuitively understand the behavior of
$\K_n^{\tilde{\P^\epsilon}}$ with Figure \ref{fig}.
In Figure \ref{fig}, the value beside a point denotes the value of
$\K_n^{\tilde{\P^\epsilon}}$ at that situation and the value beside a
diagonal line indicates the payoff Skeptic obtains in the next round.

\bigskip
As seen in (\ref{Feb 25 22:53:19 2007}), however small $\epsilon$ is,
$1+\K_n^{\tilde{\P}^\epsilon}$ will be negative if Reality makes $|s_n|$
large enough. Here we
consider the way to avoid the bankruptcy. The condition for
$\K_n^{\tilde{\P^\epsilon}}$ to be greater or equal to $-1$ is
\begin{equation}
  \K_n^{\tilde{\P^\epsilon}}=\frac{2}{\epsilon}(n-s_n^2) \geq -1
\  \Leftrightarrow \ 
|s_n|\leq\sqrt{n+\frac{2}{\epsilon}}.\label{Feb 26 16:01:43 2007}
\end{equation}
Here, suppose that Skeptic follows $\tilde{\P}^\epsilon$ and is going to announce the $n$-th move $M_n$ at the
$n$-th round. He can refer to $s_{n-1}$ but not to $s_n$. If
\[
  |s_{n-1}| > \sqrt{n+\frac{2}{\epsilon}}-1,
\]
then he should stop following $\P^\epsilon$, or else Reality can make
him bankrupt.
We let $\P^\epsilon$ denote the strategy that follows 
$\tilde{\P}^\epsilon$ under this stopping rule:
\[
\P^\epsilon : M_n = \left\{
\begin{array}{ll}
 -\epsilon s_{n-1} \quad & \mathrm{if}\
  \displaystyle|s_{i-1}|\leq\sqrt{i+\frac{2}{\epsilon}}-1, \ i=1, 2,
  \ldots, n, \\ 
 0&\mathrm{otherwise}.
\end{array}
\right.
\]

\begin{lemma}\label{lem3}
 The strategy $\P^\epsilon$ weakly forces the following $\eventthree^\epsilon$:
 \[
  \eventthree^\epsilon :=\left\{\xi\ \left|\ \exists
 n:|s_n|>\sqrt{n+1+\frac{2}{\epsilon}}-1\ \underline{\mathrm{or}}\
 \left(\limsup_{n\rightarrow\infty}|s_n|=\infty\ \&\ s_n\neq 0\ \mathrm{a.a.}\right)\right.\right\}.
 \]
\end{lemma}

\begin{proof}
 Fix any path $\xi \notin \eventthree^\epsilon$. Then
 \[
  |s_{n-1}|\leq\sqrt{n+\frac{2}{\epsilon}}-1, \quad \forall n
 \]
 holds and the capital process $\K^{\P^\epsilon}_n(\xi)$ is equal to
$\K^{\tilde{\P}^\epsilon}_n(\xi)$. Furthermore at least 
one of the following two cases holds:
\begin{enumerate}
 \item There exists $L$ such that $|s_n|<L$, $\forall n$, 
 \item $s_n=0$ holds for infinitely many $n$.
\end{enumerate}

\paragraph{case 1)}
 Since $s_n^2<L^2$, 
 \[
  \K^{\P^\epsilon}_n =
 \frac{\epsilon}{2}(n-s_n^2)>\frac{\epsilon}{2}(n-L^2)\rightarrow\infty\quad 
 (n\rightarrow \infty).
 \]
\paragraph{case 2)}
When $s_n=0$, 
 \[
  \K^{\P^\epsilon}_n =
 \frac{\epsilon}{2}n.
 \]
Therefore if $s_n=0$ occurs infinitely often, 
$\displaystyle\limsup_{n\rightarrow\infty}\K_n^{\P^\epsilon}=\infty$.
\end{proof}

Now we define $\eventthree$ by 
\[
 \eventthree :=
 \left\{\xi\ \left|\ |s_n|>\sqrt{n}-1\
 \mathrm{i.o.}\ \underline{\mathrm{or}}\
 \left(\limsup_{n\rightarrow\infty}|s_n|=\infty\ \&\ s_n\neq 0\ \mathrm{a.a.}\right)\right.\right\}.
\]
Then 
\[
 \bigcap_{i\in\mathbb{N}}\eventthree^{2^{-i}}\subset \eventthree 
\]
and by Lemma 3.2 of \cite{sv} the next corollary holds.
\begin{corollary}\label{cor1}
 Skeptic can weakly force $\eventthree$.
\end{corollary}

\subsubsection{A  strategy weakly forcing
boundedness  or two-sided unboundedness of the sum process
}
\label{sec2-onesided}

Here we consider weakly forcing the following event
\[
 \eventfour :=
 \left\{\xi\ \left|\ \limsup_{n\rightarrow \infty} |s_n| <\infty\
 \ulrmor\ \left(\limsup_{n\rightarrow\infty} s_n = \infty\ \&\
 \liminf_{n\rightarrow \infty} s_n = -\infty\right) \right.\right\}.
\]
Actually in weakly forcing $\eventfour$ we combine two one-sided strategies.
Consider the following very simple additive one-sided strategy:
\begin{equation}
\label{eq:onesided-primitive-strategy}
 \P^{-N}:M_n =
\begin{cases}
   \displaystyle\frac{1}{N} \quad & \mathrm{if} 
 \displaystyle\min_{1\leq i\leq n-1}s_i > -N, \\  
  0&\mathrm{otherwise.}
\end{cases}
\end{equation}
$\P^{-N}$ bets the constant amount $1/N$ until $s_n$ reaches $-N$
for the first time. Similarly define $\P^{+N}$ by
\begin{equation}
\label{eq:onesided-primitive-strategy-positive}
 \P^{+N}:M_n = 
\begin{cases}
   \displaystyle - \frac{1}{N} \quad & \mathrm{if} 
 \displaystyle\max_{1\leq i\leq n-1}s_i < N, \\
  0&\mathrm{otherwise.}
\end{cases}
\end{equation}
Corresponding to these strategies
in the following lemma we consider two one-sided events.
\begin{lemma}\label{lem4}
 Let $N$ be any positive number and define $\eventfour^{-N}, \eventfour^{+N}$ as follows:
\begin{eqnarray*}
 \eventfour^{-N} &:=&
 \left\{\xi\ \left|\ \min_{i=1,2,\ldots}s_i \leq -N\
 \ulrmor\ \limsup_{n\rightarrow\infty} s_n < \infty
			  \right.\right\},\\
 \eventfour^{+N} &:=&
 \left\{\xi\ \left|\ \max_{i=1,2,\ldots}s_i \geq N\
 \ulrmor\ \liminf_{n\rightarrow\infty} s_n > -\infty
			  \right.\right\}.
\end{eqnarray*}
Skeptic can weakly force $\eventfour^{+N}$ and $\eventfour^{-N}$.
\end{lemma}

\begin{proof}
Clearly the capital process $\K^{\P^{-N}}$ of $\P^{-N}$ 
is given by
 \[
  \K^{\P^{-N}}_n =
\begin{cases}
   \displaystyle\frac{s_n}{N} \quad &\mathrm{if}\ \displaystyle\min_{1\leq i\leq
    n-1}s_i > -N, \\
  -1&\mathrm{otherwise}.
\end{cases}
 \]
 This shows that $\P^{-N}$ weakly forces $\eventfour^{-N}$. The proof for
 $\eventfour^{+N}$ is almost the same by using $\P^{+N}$.
\end{proof}

\begin{corollary}\label{cor2}
Skeptic can weakly force $\eventfour$.
\end{corollary}

\begin{proof}
Define $\eventfour^-$ and $\eventfour^+$ as follows:
\begin{eqnarray*}
  \eventfour^- &:=&
 \left\{\xi\ \left|\ \liminf_{n\rightarrow \infty} s_n =-\infty\
 \ulrmor\ \limsup_{n\rightarrow\infty} s_n < \infty \right.\right\},\\
   \eventfour^+ &:=&
 \left\{\xi\ \left|\ \liminf_{n\rightarrow \infty} s_n >-\infty\
 \ulrmor\ \limsup_{n\rightarrow\infty} s_n = \infty \right.\right\}.\\
\end{eqnarray*}
Then we can write 
 \begin{eqnarray*}
  \eventfour^- &=& \bigcap_{N\in\mathbb{N}}\eventfour^{-N}, \\
  \eventfour^+ &=&\bigcap_{N\in\mathbb{N}}\eventfour^{N},
\end{eqnarray*}
and 
 \[
  \eventfour = \eventfour^-\cap \eventfour^+.
 \]
 By the Lemma 3.2 of \cite{sv}, the corollary holds.
\end{proof}

\subsubsection{Proof of Theorem \ref{thm2}}

Using Corollary \ref{cor1} and Corollary \ref{cor2}, now we can
prove Theorem \ref{thm2}.

\medskip
\begin{proof}
 From Corollary \ref{cor1} and Corollary \ref{cor2}, Skeptic can weakly
 force $\eventthree\cap \eventfour$. So we only have to show $\eventtwo\supset \eventthree\cap \eventfour$.
 Now we set
 \begin{eqnarray*}
  A_1 &:=& \left\{\xi\ \left|\ \limsup_{n\rightarrow\infty}|s_n|=\infty\ \&\
		      s_n\neq 0\ \mathrm{a.a.}\right.\right\},\\
  A_2 &:=& \left\{\xi\ \left|\ \limsup_{n\rightarrow \infty} |s_n|
		      <\infty\right.\right\},\\
  A_3 &:=& \left\{\xi\ \left|\ \limsup_{n\rightarrow\infty} s_n = \infty\ \&\
 \liminf_{n\rightarrow \infty} s_n = -\infty \right.\right\}.
 \end{eqnarray*}
 Then we can write
\[
  \eventthree = \eventtwo\cup A_1,\qquad 
  \eventfour = A_2\cup A_3.
\]
By definition
\[
\emptyset =  \eventtwo\cap A_2 =   A_1\cap A_2 
=  A_1\cap A_3 
\]
and therefore
 \begin{eqnarray*}
  \eventthree\cap \eventfour = \eventtwo\cap A_3\subset \eventtwo.
 \end{eqnarray*}
\end{proof}

This proof also shows the next theorem.
\begin{theorem}\label{thm3}
 Skeptic can weakly force $A_3$.
\end{theorem}

\section{One sided strategies}
\label{sec:one-sided}
The statement of Theorem \ref{thm2} 
is only concerned with the behavior of $|s_n|$
and Reality is forced
to make $|s_n| > \sqrt{n}-1$ infinitely often. But it says
nothing about the sign of $s_n$. Hence Reality can choose a path such
that $s_n > \sqrt{n}-1$ infinitely often but $s_n <-\sqrt{n}+1$ only
finitely often. In this section we prove the following theorem which eliminates
this shortcoming.
 \begin{theorem}\label{thm4}
  Skeptic can weakly force the following $\eventfive$ and $\eventsix$:
  \begin{eqnarray*}
   \eventfive &:=& \{\xi \mid s_n>\sqrt{n}-1\ \mathrm{i.o.}\},\\
   \eventsix &:=& \{\xi \mid s_n<-\sqrt{n}+1\ \mathrm{i.o.}\}.\\
  \end{eqnarray*}
 \end{theorem}

The statement in Theorem  \ref{thm4} seems to be innocuous and
one might expect that it can be proved by the obvious symmetry
of the fair-coin game.  Actually we found it difficult
to prove Theorem  \ref{thm4} by combination of simple strategies.
Recall that in the previous section, except for combining countably
many strategies, the individual strategies were very simple and
explicit.  Furthermore it should be possible to  generalize 
the results in the previous section  
to more general protocols than the fair-coin game
by introducing pricing of quadratic  hedges as in Chapter 4 of
\cite{sv}.

On the other hand, our proof of Theorem \ref{thm4} uses the fact that
in the fair-coin game it is conceptually very easy to determine the
price of every variable.  Mathematically it is the same as the pricing
of options for binomial models, which is explained in standard
introductory textbooks on mathematical finance (e.g.\ \cite{finance}).
See also \cite{ts} for a game-theoretic exposition of the pricing
formulas for the binomial model.

\subsection{Two stopping times}
\begin{figure}[h]
\begin{center}
  \includegraphics{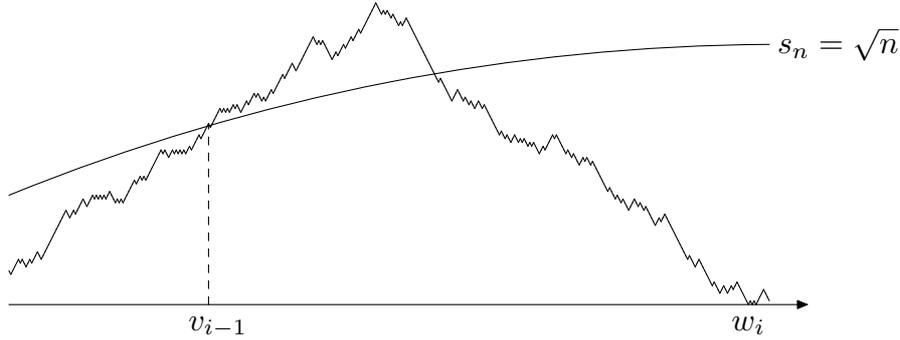}
\end{center}
 \caption{definition of $w_i$ and $v_i$}\label{fig2}
\end{figure}

For $i = 1, 2, \ldots$ , we define stopping times $w_i$ and $v_i$ by
\begin{eqnarray*}
 w_i &:=& \min\left\{n>v_{i-1}\ \left|\ s_n=0\right.\right\},\\
 v_i &:=& \min\left\{n>w_i\ \left|\ |s_n|>\sqrt{n}-1\right.\right\},
\end{eqnarray*}
where $v_0:=0$ and if the set in the definition is empty then the
value of the variable is $\infty$. 
$v_i$ is the first hitting time of the two-sided $\sqrt{n}$-boundary 
after leaving the origin at $w_i$ and $w_i$ is the first time of
returning to the origin after $v_{i-1}$. See Figure \ref{fig2}.

Now for  $i=1,2,\ldots$ , we define a variable $X_i$ by
\[
 X_i := \left\{
 \begin{array}{ll}
  1&v_i<\infty\ \&\ s_{v_i}<0, \\
  0&\mathrm{otherwise}.
 \end{array}
 \right. 
\]
We can think of $X_i$ as a ticket, which pays you one dollar if the
sum process hits the negative boundary $-\sqrt{n}$ (rather than the
positive boundary $\sqrt{n}$\,) at time $v_i$.

\subsection{Proof of theorem \ref{thm4}}
This section is devoted to the proof of Theorem \ref{thm4}.
It is conceptually very easy.
The essential point is the proof of Lemma \ref{lem5} below. We begin
by giving our proof other than Lemma \ref{lem5}.

Remember that Reality is forced the event $A_3\cap \eventtwo$. 
Therefore we can assume that 
she chooses $\xi\in A_3\cap \eventtwo$.  Therefore for proving
Theorem \ref{thm4}, it suffices to prove 
that Skeptic can weakly force $A_3\cap \eventtwo\Rightarrow
\eventfive$ (cf.\ Lemma 2.1 of \cite{ktt-unbounded}). 

In Lemma \ref{lem5}, we prove $\overline{\mathrm{E}}_t[X_i]
\leq \frac{1}{2}$ for any situation $t\in [w_i]$ by constructing the
replicating strategy of $X_i$. We let $\P^{X_i}$ denote this
strategy. 
Once $\P^{X_i}$ is constructed, the strategy weakly forcing 
$A_3\cap \eventtwo
\Rightarrow \eventfive$ is given as follows:
\begin{itemize}
 \item Buy $\frac{1}{2}\K X_i$ when the present situation is in $[w_i]$
       for $i=1, 2, \ldots$,
\end{itemize}
where $\K$ denotes the present capital Skeptic possesses and buying
$\frac{1}{2}\K X_i$ actually means running $\frac{1}{2}\K\P^{X_i}$.

Suppose Reality chooses the path $\xi\notin (A_3\cap \eventtwo)^C\cup
\eventfive$. Since $\xi\in A_3\cap \eventtwo$, $w_i<\infty$ and
$v_i<\infty$ for any
$i$. Thus Skeptic runs $\frac{1}{2}\K_{w_i}\P^{X_i}$ from the $w_i$-th round
for each $i$. After the $v_i$-th round, 
his capital becomes
$\frac{3}{2}\K_{w_i}$ if $s_{v_i}<0$ and $\frac{1}{2}\K_{w_i}$ if
$s_{v_i}>0$. But $s_{v_i}<0$ for all sufficiently large $i$ since
$\xi\notin E_5$, then Skeptic's capital
increases to $\infty$.

Now it remains to prove the following lemma.

\begin{lemma}\label{lem5}
 For any $i \in \mathbb{N}$ and any $t\in [w_i]$
 \begin{equation}
  \overline{\mathrm{E}}_t[X_i]
\leq \frac{1}{2}.\label{Mar 21 00:03:49 2007}
 \end{equation}
\end{lemma}

\begin{proof}
 First we rephrase the lemma for simplicity.
 Fix any situation $t\in[w_i]$ and suppose that the length of $t$ is $l$.
 Define a stopping time $u_l$ and
 a variable $Y_l$ as
 \begin{eqnarray*}
  u_l &:=& \min\left\{n\ \left|\ |s_n|>\sqrt{n+l}-1\right.\right\},\\
  Y_l &:=& \left\{
 \begin{array}{ll}
  1&\mathrm{if}\ u_l<\infty\ \&\ s_{u_l}<0, \\
  0&\mathrm{otherwise}.
 \end{array}
 \right.
 \end{eqnarray*}
 Considering the upper price of $X_i$ at the situation $t$ is equivalent to
 considering the upper price of $Y_l$ at the situation $\Box$, since
 $s(t)=0$. Thus it suffices to show $\overline{\mathrm{E}}[Y_l] \leq
 \frac{1}{2}$ for the proof of
 $\overline{\mathrm{E}}_t[X_i]\leq\frac{1}{2}$.

\bigskip
The set $\{Y_l=1\}$ can be decomposed as
\begin{eqnarray*}
 \{\xi \mid Y_l(\xi)=1\}&=&\bigcup_{i=1}^\infty A_i,\\
 A_i &:=& \{\xi \mid u_l(\xi)=i, s_i<0\}.
\end{eqnarray*}
Using $A_i$, $Y_l$ can be decomposed as
\[
 Y_l = \sum_{i=1}^{\infty}Z_i ,  \qquad 
 Z_i(\xi) := \left\{
\begin{array}{ll}
 1& \mathrm{if}\ \xi\in A_i, \\
 0&\mathrm{otherwise}.
\end{array}
\right.
\]

Whether $\xi\in A_i$ or not depends solely on $\xi_i$, so
 $A_i$ can be decomposed into the cylinder sets defined by the
 situations of length $i$, that is, there exist $t^i_j\ (j=1, 2,
 \ldots, a_i)$ such that
\[
 A_i = \bigcup_{j=1}^{a_i}O_{t^i_j},
\]
where the length of $t^i_j$ is $i$\ ($j=1, 2, \ldots, a_i$).
We set, 
\[
 t^i_j = y^{i, j}_1y^{i, j}_2\cdots y^{i, j}_i,\quad  (y^{i, j}_p \in \{-1,
 1\}, p = 1, 2, \ldots, i).
\]
Using $t^i_j$, we define the strategy $\P^{i, j}$ by
\[
 \P^{i, j}:\ M_n =\left\{
 \begin{array}{ll}
  y^{i, j}_n\K_{n-1}& \mathrm{for}\ n=1, 2, \ldots, i,\\
  0& \mathrm{for}\ n>i,
 \end{array}
 \right.
  \]
where we temporarily suppose $\K_{0}=2^{-i}$. Intuitively speaking, this strategy prepares the amount 
of $2^{-i}$ as the
 initial capital and bet all the available capital on the realization of the
 situation $t^i_j$. Hence, if $t^i_j$ realizes, that is, $x_p = y^{i, j}_p,\ (p=1, 2, \ldots, i)$, then the capital grows to
 $2^{i}$ times, otherwise the capital becomes zero. Thus the capital process
 $\K^{\P^{i, j}}$ satisfies
\[
 2^{-i}+\K_i^{\P^{i,j}}(\xi)=\left\{
 \begin{array}{ll}
  1&\mathrm{if}\ \xi\in O_{t^i_j}, \\
  0&\mathrm{otherwise}.
 \end{array}
 \right.
\]
Moreover, we define the strategy $\P^i$ by
 \[
  \P^i := \sum_{j=1}^{a_i}\P^{i, j}.
 \]
 The strategy $\P^i$ requires the amount of 
$a_i2^{-i}$ as the initial capital and its capital process is written as
 \[
 a_i2^{-i} + \K_i^{\P^i}(\xi) =\left\{
 \begin{array}{ll}
  1& \mathrm{if}\ \xi\in \bigcup_{j=1}^{a_i}O_{t^i_j}, \\
  0&\mathrm{otherwise}.
 \end{array}
 \right.
 \]
Thus Skeptic can replicate $Z_i$ with the initial capital $a_i2^{-i}$.
\medskip
Then, Skeptic can replicate $\displaystyle Y_l=\sum_{i=1}^\infty Z_i$
 with initial capital
 $\displaystyle\sum_{i=1}^{\infty}a_i2^{-i}$, so it suffices to show
\[
  \sum_{i=1}^{\infty}a_i2^{-i}\leq \frac{1}{2}.
\]
 for the
 proof of $\overline{\mathrm{E}}[Y_l]\leq\frac{1}{2}$
 
\bigskip
Fix an arbitrary large number $k$. We consider the event $B_k$ defined by
\begin{eqnarray*}
 B_k &:=& \{\xi \mid u_l(\xi)\leq k, s_l < 0\}.
\end{eqnarray*}
Whether $\xi\in B_k$ or not depends solely on $\xi_k$, so $B_k$ can be
 decomposed into the cylinder sets:
 \begin{eqnarray}
 B_k = \bigcup_{q=1}^{b_k}O_{t'^k_q},\label{Jan 25 17:48:57 2007}
 \end{eqnarray}
 where the length of $t'^k_q$ is $k$\ ($q=1, 2, \ldots, b_k$).
 Here we show $b_k\leq2^{k-1}$. First, remember that there are just
 $2^k$ situations of length $k$. If we define $B_k^-$ as
 \[
  B_k^- := \bigcup_{q=1}^{b_k}O_{-t'^k_q},
 \]
 then $B_k\cap B_k^-=\emptyset$ by definition. Thus $b_k\leq
 \frac{1}{2}\cdot2^k=2^{k-1}$. 
 Furthermore $B_k$ can be decomposed also into $A_i$:
\[
 B_k = \bigcup_{i=1}^{k}A_i=\bigcup_{i=1}^{k}\bigcup_{j=1}^{a_i}O_{t^i_j}.
\]

Just $2^{k-i}$ situations of $\{t'^k_q\}_{q=1}^{b_k}$ follow $t^i_j$,
 so the cylinder set $O_{t^i_j}$ can be decomposed into $2^{k-i}$
 cylinder sets:
 \[
 O_{t^i_j}=\bigcup_{h=1}^{2^{k-i}}O_{t'^k_{c[i,j, h]}},
\]
 where $1\leq c[i, j, h]\leq b_k$ and $c[i,j, h]\neq c[i', j', h']$ if $(i, j, h)\neq (i', j', h')$.
Thus, 
\begin{equation}
 B_k =
  \bigcup_{i=1}^{k}\bigcup_{j=1}^{a_i}\bigcup_{h=1}^{2^{k-i}}O_{t'^k_{c[i,j,
  h]}}.\label{Jan 25 17:49:07 2007}
\end{equation}
By (\ref{Jan 25 17:48:57 2007}) and (\ref{Jan 25 17:49:07 2007})
\[
 b_k=\sum_{i=1}^{k}a_i 2^{k-i}.
\]
Since $b_k\leq 2^{k-1}$, 
\begin{eqnarray*}
 \sum_{i=1}^{k}a_i 2^{k-i}\leq 2^{k-1} 
\quad  \Rightarrow \quad 
\sum_{i=1}^{k}a_i 2^{-i}\leq\frac{1}{2}
\quad  \Rightarrow \quad
\sum_{i=1}^{\infty}a_i 2^{-i}\leq\frac{1}{2}.
\end{eqnarray*}
\end{proof}

Lastly, we show that the inequality in the (\ref{Mar 21 00:03:49
2007}) is in fact an equality. Here we decompose $\Omega$ into three subsets,
\begin{eqnarray*}
 D^i_1 &:=& \{\xi\mid v_i<\infty\ \&\ s_{v_i}<0\},\\
 D^i_2 &:=& \{\xi\mid v_i=\infty\},\\
 D^i_3 &:=& \{\xi\mid v_i<\infty\ \&\ s_{v_i}>0\}.\\
\end{eqnarray*}
Then $X_i=1_{D^i_1}$, where $1_{D^i_1}$ is the indicator function of
$D^i_1$. By the definition of the  upper price,
\begin{equation}
  1=\overline{\mathrm{E}}_t[1]\leq\overline{\mathrm{E}}_t[1_{D^i_1}]+\overline{\mathrm{E}}_t[1_{D^i_2}]+\overline{\mathrm{E}}_t[1_{D^i_3}]\label{Mar 23 03:04:53 2007}
\end{equation}
By symmetry property,
\begin{equation}
  \overline{\mathrm{E}}_t[1_{D^i_1}]=\overline{\mathrm{E}}_t[1_{D^i_1}].\label{Mar 23 03:05:02 2007} 
\end{equation}
Since Skeptic can force $v_i<\infty$ we have
\begin{equation}
  \overline{\mathrm{E}}_t[1_{D^i_2}] = 0.\label{Mar 23 03:05:17 2007}
\end{equation}
Equations (\ref{Mar 23 03:04:53 2007}), (\ref{Mar 23 03:05:02 2007}) and
(\ref{Mar 23 03:05:17 2007}) shows 
$\overline{\mathrm{E}}_t[1_{D^i_1}]\geq\frac{1}{2}$, thus $\overline{\mathrm{E}}_t[1_{D^i_1}]=\overline{\mathrm{E}}_t[X_i]=\frac{1}{2}$.

\section{Some discussions}
\label{sec:discussion}

In this paper we showed that Skeptic can (weakly) force $\eventone,
\eventtwo, \eventfive$ and $\eventsix$ in the fair-coin game. 
As mentioned in Section \ref{sec:intro}
these statements are weaker than LIL,  which is
shown in \cite{sv} in the game-theoretic framework. But we want to
emphasize the simplicity of our strategies. Actually, Skeptic needs only
to keep the value of $s_n$ in memory in the strategies forcing $\eventone$ and
$\eventtwo$. 

\smallskip
In the proof of Lemma \ref{lem5}, we only proved the existence of the
replicating strategy of $Y_l$ rather than providing an  explicit
formula for  bet (the move $M_n$ of Skeptic) of the strategy. 
The bet of the replicating strategy is
directly related with the price of $Y_l$ at an arbitrary situation by
the argument of ``delta hedge''\cite{ts}. Let $\eta(n, s)$ denote the price
of $Y_l$ given the round $n$ and the value of process $s$. Here let us
consider the problem in the measure-theoretic framework rather than the 
game-theoretic framework. 
Then $\eta(n, s)$ can
be written as the measure-theoretic conditional expectation:
\begin{equation}
\label{eq:conditional-Y}
\eta(n, s) = \mathrm{E}[Y_l \mid s_n = s].
\end{equation}
Given $\eta(n, s)$, the bet of replicating strategy by
delta hedge is calculated as follows(\cite{ts}):
\[
 M_n = \frac{\eta(n+1, s+1)-\eta(n+1, s-1)}{2}.
\]
But in practice it would be difficult to express $\eta(n, s)$
analytically.  For the case of Brownian motion \cite{novikov}
gives results on (\ref{eq:conditional-Y}).  However they are very
complicated involving zeros of a special function.

In order to prove the existence of the replicating strategy, we used the
argument of betting on specific paths. This type of argument can be 
found in the field of algorithmic theory of randomness, for
instance Muchnik et al.\ uses the same idea in \cite[Theorem 9.4]{muchnik}. 
We think that the idea is logically very powerful because it can be used
to prove the existence of a superreplicating strategy for any ticket in
the fair-coin game.

\end{document}